\documentclass{amsart}
\begin{document}
\title{Biography of John Rainwater}
\author{Robert R. Phelps}
\address{Department of Mathematics, University of Washington}
\thanks{This essay is published in \emph{Topological Commentary}
\textbf{7} no.\ 2 (2002).}
\thanks{This is adapted from the introduction to a volume of J.R.'s
collected works [\oldstylenums{18}], on file at the University of
Washington.}
\date{October, \oldstylenums{2002}}
\maketitle

The following paragraphs will describe the origins of John Rainwater, the
impact of his work, the motivations for various parts of it and the
prospects for his future.

\section*{Origins}

John Rainwater came into existence at the University of Washington in
\oldstylenums{1952} when Nick Massey, a mathematics graduate student in
Prof.\ Maynard Arsove's beginning real variables class, erroneously
received a blank registration card.  (In those years, each student filled
out a card for every class, which first circulated among various
tabulating clerks in the registrar's office before being sent to the
professor.)  He and a fellow graduate student, Sam Saunders, decided to
use the card to enroll a fictional student, and since it was raining at
the time, decided to call him ``John Rainwater''.  They handed in John
Rainwater's homework regularly, so it wasn't until after the first midterm
exam that Prof. Arsove became aware of the deception.  He took it well,
even when he later opened an ``exploding'' fountain pen with John
Rainwater's name engraved on it which had been left on the classroom
table.  After remarks by Arsove, such as ``I guess I'll never see
Rainwater except in a barrel,'' virtually all the students learned of the
Rainwater prank.  Several years later, when a group of mathematics
graduate students and junior faculty solved some American Mathematical
Monthly problems over coffee in the student union, it was only natural
that they would submit them to the Monthly under J.R.'s name.  As it does
at present, the {\scshape maa} doubtless suggested that someone who was
submitting problem solutions should become a member of the Association.
At that time, a prospective \textsc{maa} member had to obtain the
endorsement of two current members.  An ideal endorser would have been
Professor Carl Allendoerfer, who was not only Chairman of the Mathematics
Department but also President of the \textsc{maa}.  Unfortunately, he was
not the kind of person to participate in a joke involving the
\textsc{maa}.  The solution was simple:  A secretary was persuaded to put
the endorsement form in with a pile of papers which required the busy
Chairman's signature, and the deed was accomplished.  (Allendoerfer was
reported to have been more than a bit unhappy when he eventually learned
of the matter.)  John has continued to publish problem solutions; see the
paragraph on Authorship and Motivation below.

The first of John Rainwater's ten published research papers were written
in \oldstylenums{1958} and \oldstylenums{1959} by John Isbell, a young
Assistant Professor.  Isbell's response to queries concerning his
motivation for using J.R.\ as a pseudonym has been simply to quote
Friedrich Schiller ``Der Mensch ist nur da ganz Mensch, wo er spielt.''
Unlike many of the Rainwater papers, there is no hint within these first
two papers as to the actual authorship. (For present-day readers of this
account, it should be pointed out that academic jobs were plentiful in the
late \oldstylenums{50}'s, so amassing a publication record was not as
pressing as it is today.)  A functional analysis seminar was started in
the late \oldstylenums{1960}'s, and it seemed appropriate to call it the
``John Rainwater Seminar'', in view of the work that J.R.\ had done in
that field.  As the functional analysts in the department over the years
have changed fields, died or retired, the John Rainwater seminar has
changed its emphasis from simply functional analysis to functional
analysis plus Fourier analysis, then to these two fields plus dynamical
systems.  At present it is primarily a seminar in the latter area.

\section*{Impact}

What about the impact of John Rainwater's research?  This is something
that it is always difficult to measure, but a rough idea can be obtained
from the citations listed in the Institute for Scientific Information's
``Web of Science''.  (The numbers which follow are those in the
bibliography below.) His most cited paper, with \oldstylenums{19}
citations, the latest in \oldstylenums{2000}, was the first one, in
topology.  Number \oldstylenums{6}, his \oldstylenums{1969} \textsc{pams}
note on Day's norm, is his next most cited paper, with \oldstylenums{17}
citations, the latest in \oldstylenums{1998}.  The algebra paper (number
\oldstylenums{14}) has \oldstylenums{14} citations, the latest in
\oldstylenums{2001}.  There are eight citations to both his second paper,
the latest in \oldstylenums{1989}, and his third one, the one--page
\oldstylenums{1963} \textsc{pams} note, the latest in \oldstylenums{2002}.
The latter has also been cited in books as ``Rainwater's Theorem''.
Paper \oldstylenums{15} on convex functions has six citations, the latest
in \oldstylenums{1995}, while paper \oldstylenums{10} on regular matrices
has four citations, the latest in \oldstylenums{1997}. Paper
\oldstylenums{7}, also about Day's norm, has three citations, the latest
in \oldstylenums{1993} and number \oldstylenums{8} has two, the latest in
\oldstylenums{1998}. There is even one citation to number
\oldstylenums{13}, his unpublished \oldstylenums{1967} Rainwater Seminar
note on Lindenstrauss spaces. Papers \oldstylenums{4} and
\oldstylenums{16} have received no citations to date.  In summary, it
appears that most of John Rainwater's published work has been reasonably
well received.

\section*{Authorship and motivation}

As mentioned earlier, John Isbell wrote the first two John Rainwater
papers.  I wrote the third one, and the motivation was simple:  I needed
the theorem for a paper I was writing, but it was a ``folk theorem'',  
observed more or less independently by five other mathematicians.  It  
would have been a bit silly to have a half-page paper appear under the
names of five authors, so I wrote to each of them, obtaining their
approval to publish it under J.R.'s name (who acknowledged their
``extremely useful conversations'').  The fourth paper was written by
Irving Glicksberg.  I'm not certain of his motivation; this might be the
only instance where the actual author thought that the paper was not good
enough to appear under his own name.  (The absence of any citations for it
suggests that he was right.)  Paper \oldstylenums{5} is an unpublished
note that was circulated to the experts.  Paper \oldstylenums{6} about
Day's norm was again a case of a large number of people (six this time)
having provided various parts of the proof, so a pseudonym was very much 
in order.  Again, they are thanked for ``helpful discussions''. Paper  
\oldstylenums{7}, which deduced additional properties of Day's norm, was
written by Edgar Asplund for a conference in Aarhus.  Paper
\oldstylenums{8} on the abstract F. and M. Riesz theorem was written by
Glicksberg.  The result was proved by four people but it was only one and
a half pages long, so the motivation is again obvious.  The unpublished  
seminar note in number \oldstylenums{9} was written by me for the
Rainwater Seminar. Number \oldstylenums{10} (\emph{Regular matrices
\ldots}) was written by John Isbell again (who, incidentally, has authored
or co-authored six other pseudononymous papers under two other names).  I
proved the result in number \oldstylenums{11} and subsequently learned  
that Day and Pelczynski were also aware of it.  John Giles also proved it
independently, and he was rather unhappy that the rest of us thought that
it was not important enough to publish.  Number \oldstylenums{12} is an
exposition I wrote for the Rainwater Seminar.  Number \oldstylenums{13} 
was proved by Peter D. Morris and me when he was visiting here; we thought
it was too clearly an example of ``proof by theorem'' to be worth
publishing.  Paper \oldstylenums{14} is a departure for John Rainwater.
Not only is it in algebra, but he doesn't thank anyone for helpful
conversations.  He notes, however, that his work was supported by four  
different grants.  (Culprits this time were Ken Brown, Ken Goodearl, Toby
Stafford and Bob Warfield.) Number \oldstylenums{15} was written by Isaac 
Namioka and me.  We were able to generalize Elena Verona's theorem simply
because we were familiar with some techniques she (as a new Ph.D.)  had
not yet learned and we would have been embarrassed to take personal credit
for that.  Number \oldstylenums{16} was written by me and David Preiss; it
has not yet shown up in in the Science Citation Index.  Number
\oldstylenums{17} is a collection (surely incomplete) of problems or
solutions which J.R.\ has published in the American Mathematical Monthly,
the earliest in \oldstylenums{1959} (sent by John Isbell), the latest in
\oldstylenums{1994} (sent by me).  I have no idea who was responsible for
problem \oldstylenums{4908} (\oldstylenums{1963}) from the University of
British Columbia.

\section*{Future}

Where does John Rainwater go from here?  It would be a shame if he were to
die.  He is not as old or famous as N. Bourbaki (who may still be alive)
but he is clearly older than Peter Orno, who only has three publications
to his name, all in the \oldstylenums{1970}'s.  (At least one of his
authors had an interest in pornography, hence P. Orno.)  He is also older
than M. G. Stanley (with four papers) and H.~C.~Enoses (with only two). 
It is to be hoped that someone will be able to help John Rainwater carry
on, so that in the future people won't ask ``Who killed J.R.?'' The
previous paragraph gives some guidance to potential Rainwater authors.   
Clearly, in the present employment climate, no junior faculty member would
be willing to emulate John Isbell and publish a really good paper under a 
pseudonym. The other successful J.R.\ papers were those which simply had
too many authors, all of whom were securely tenured faculty, so keep him 
in mind if you find yourself in such a situation.  It would be nice to   
keep the tradition going.

\subsection*{Acknowledgements}

A number of people have helped me to write this biography.  My
appreciative thanks go to Maynard Arsove, John Isbell and Sam Saunders for
their information on J.R.'s origins and to the University of
Washington Mathematics Librarian Martha Tucker for furnishing the citation
material.

\section*{Bibliography}

\renewcommand{\labelenumi}{\oldstylenums{\theenumi}.}


\begin{enumerate}
\setlength{\itemindent}{-0.3in}

\item
\emph{Spaces whose finest uniformity is metric},
Pacific J.~Math.\
\textbf{\oldstylenums{9}} (\oldstylenums{1959}),
\oldstylenums{567}--\oldstylenums{570}.
\textsc{mr}~\oldstylenums{21}\#\oldstylenums{5180};
\textsc{zbl}~\oldstylenums{0088}.\oldstylenums{38301}

\item
\emph{A note on projective resolutions},
Proc.~Amer.~Math.~Soc.\
\textbf{\oldstylenums{10}} (\oldstylenums{1959}),
\oldstylenums{734}--\oldstylenums{735}.
\textsc{mr}~\oldstylenums{23}\#A\oldstylenums{618};
\textsc{zbl}~\oldstylenums{0105}.\oldstylenums{16403}

\item
\emph{Weak convergence of bounded sequences},
Proc.~Amer.~Math.~Soc.\
\textbf{\oldstylenums{14}} (\oldstylenums{1963}),
\oldstylenums{999}.
\textsc{mr}~\oldstylenums{27}\#\oldstylenums{5111};
\textsc{zbl}~\oldstylenums{0117}.\oldstylenums{08302}

\item
\emph{A remark on regular {B}anach algebras},
Proc.~Amer.~Math.~Soc.\
\textbf{\oldstylenums{18}} (\oldstylenums{1967}),
\oldstylenums{255}--\oldstylenums{256}.
\textsc{mr}~\oldstylenums{34}\#\oldstylenums{8223};
\textsc{zbl}~\oldstylenums{0171}.\oldstylenums{33802}

\item
\emph{On a renorming theorem of {K}lee},
Unpublished note,
\oldstylenums{1968}.

\item
\emph{Local uniform convexity of {D}ay's norm on
{$c_{\oldstylenums{0}}(\Gamma)$}},
Proc.~Amer.~Math.~Soc.\
\textbf{\oldstylenums{22}} (\oldstylenums{1969}),
\oldstylenums{335}--\oldstylenums{339}.
\textsc{mr}~\oldstylenums{39}\#\oldstylenums{4647};
\textsc{zbl}~\oldstylenums{0185}.\oldstylenums{37602}

\item
\emph{{D}ay's norm on
{$c_{\oldstylenums{0}}(\Gamma)$}},
Proc. of the Functional Analysis Week, Aarhus  \
\textbf{\oldstylenums{8}}
 (\oldstylenums{1969}),
\oldstylenums{46}--\oldstylenums{50}, Matematisk Inst.,
Aarhus Univ., Aarhus.
\textsc{mr}~\oldstylenums{40}\#\oldstylenums{7778};
\textsc{zbl}~\oldstylenums{0235}.\oldstylenums{46048}

\item
\emph{A note on the preceding paper},
Duke Math.~J.\
\textbf{\oldstylenums{36}} (\oldstylenums{1969}),
\oldstylenums{799}--\oldstylenums{800}.
\textsc{mr}~\oldstylenums{44}\#\oldstylenums{7299};
\textsc{zbl}~\oldstylenums{0201}.\oldstylenums{45801}

\item
\emph{A characterization of certain dual unit balls},
Rainwater Sem.~Notes,  
\oldstylenums{1970}.

\item
\emph{Regular matrices with nowhere dense support},
Proc.~Amer.~Math.~Soc.\
\textbf{\oldstylenums{29}} (\oldstylenums{1971}),
\oldstylenums{361}.
\textsc{mr}~\oldstylenums{43}\#\oldstylenums{5330};
\textsc{zbl}~\oldstylenums{0213}.\oldstylenums{08301}

\item
\emph{A non-reflexive Banach space has non-smooth third conjugate space},
Rainwater Sem.\ Notes,
\oldstylenums{1972}.

\item
\emph{A theorem of Ekeland and Lebourg on Frechet differentiability of
convex functions on Banach Spaces},
Rainwater Sem.\ Notes, 
\oldstylenums{1976}.

\item
\emph{Lindenstrauss spaces which are Asplund spaces},
Rainwater Sem.\ Notes,
\oldstylenums{1976}--\oldstylenums{77}.

\item
\emph{Global dimension of fully bounded {N}oetherian rings},
Comm.~Algebra
\textbf{\oldstylenums{15}} (\oldstylenums{1987}), no.~\oldstylenums{10},
\oldstylenums{2143}--\oldstylenums{2156}.
\textsc{mr}~\oldstylenums{89}b:\oldstylenums{16032};
\textsc{zbl}~\oldstylenums{0628}.\oldstylenums{16010}

\item
\emph{Yet more on the differentiability of convex functions},
Proc.~Amer.~Math.~Soc.\
\textbf{\oldstylenums{103}} (\oldstylenums{1988}), 
no.~\oldstylenums{3}, \oldstylenums{773}--\oldstylenums{778}.
\textsc{mr}~\oldstylenums{89}m:\oldstylenums{46081};
\textsc{zbl}~\oldstylenums{0661}.\oldstylenums{49007}

\item
\emph{A class of null sets associated with convex functions on {B}anach
spaces},
Bull.~Austral.~Math.~Soc.\
\textbf{\oldstylenums{42}} (\oldstylenums{1990}), no.~\oldstylenums{2},
\oldstylenums{315}--\oldstylenums{322}.
\textsc{mr}~\oldstylenums{91}j:\oldstylenums{46050};
\textsc{zbl}~\oldstylenums{0724}.\oldstylenums{46017}

\item
\emph{Problems/Solutions published by John Rainwater}.

\item
\emph{Collected Works of John Rainwater},
Department of Mathematics, University of Washington.

\end{enumerate}
 
\end{document}